\documentclass[a4paper]{article}
\title{On the counterexamples to Borsuk's conjecture by Kahn and Kalai}
\author{Thomas Jenrich}
\date{2018-10-11}
\begin{document}
\maketitle

\section{Abstract}

In the concluding remarks of their 1993 published and now famous paper
\cite{K}, Jeff Kahn and Gil Kalai wrote in particular:

``Our construction shows that Borsuk's conjecture is false for $d = 1,325$
 and for every $d > 2,014$.''

But, as Bernulf Wei{\ss}bach remarked in his paper \cite{W} from 2000,
a simple (few steps for an upper-class pocket calculator) computation
indicates that that claim is not true for $d=1325$.

William Kretschmer (Univ. of Texas) sent me his writeup on that paper by Kahn
and Kalai, in particular pointing out that the derivation of the formula used
in that computation disregarded a certain aspect, that way missed the chance
to remove one final halving from that formula and to indeed provide a proof
that a certain point set is a counterexample for dimension 1325
(and all higher dimensions up to 1560, too).

This updated article takes a closer look at that derivation, gives an own,
much more detailed and formal version of it that delivers the
improved/corrected formula, and contains some further conclusions.

\section{Borsuk's conjecture}

In \cite{B} (1933) Karol Borsuk asked whether each bounded set in the
$n$-dimensional Euclidean space can be divided into $n$+1 parts of smaller diameter.
The diameter of a set is defined as the supremum (least upper bound) of the
distances of contained points. Implicitly, the whole set is assumed to
contain at least two points. The hypothesis that the answer to that question
is positive became famous under the name \emph{Borsuk's conjecture}.

\subsection{The special case of (0,1)-vectors with constant weight}

In \cite{L} (1984) D. Larman gave a translation of Borsuk's conjecture for
this case into the terminology of finite sets:

Let $K$ be a family of $k$-subsets of $\{1,2,...,n\}$ such that every two
members of $K$ have $t$ elements in common. Then $K$ can be partitioned into
$n$ parts so that in each part every two members have $(t+1)$ elements in
common.

\medskip

Disproving this conjecture means disproving Borsuk's conjecture, too.

\section{The counterexamples by Kahn and Kalai}

\subsection{The Frankl-Wilson theorem}

Theorem 2 (quoted below) plays an important role in \cite{K}. It is not
exactly the theorem proved in \cite{F} by P. Frankl and R. Wilson but derived
from it in few steps, mainly to let finite sets be the objects.

\medskip

Let $k$ be a prime power and $n = 4k$. Let $K$ be a family of $n/2$-subsets of $\{1,2,...,n\}$, so that no two sets
in the family have intersection of size $n/4$.
Then $ |K| \leq 2 \cdot {{n-1}\choose{n/4-1}}$ \ .

\subsection{The construction}

\subsubsection{The original description}

While the respective section of \cite{K} is clear enough with respect to
construction, cardinality and dimension of the relevant point set $K$,
the derivation of an upper bound for the cardinality of subsets of smaller
diameter is rather sketchy. It seems that not taking in account that each of
the partitions used to construct $K$ exclusively contains two subsets of $V$
(to count in the application of the Frankl-Wilson theorem) caused that the
resulting upper bound is exactly the double of what proper reasoning would
give. Consequently, the lower bound of the number of subsets of smaller
diameter is unnecessarily halved, in particular too small for the else
possible direct proof that the 1325-dimensional point set generated in the
case $k=13$ is a counterexample to Borsuk's conjecture.

\subsubsection{A more detailed and formal description}

There are propositions (definitions, premisses and conclusions) and
comments. If the symbol * starts the visible part of a line, the remaining
part of the line is a comment, with respect to the preceding line.

The strings dist, diam and dim denote the functions distance, diameter
and dimension, resp.

\textbf{Premiss} (k)

 $k$ is a prime power

\textbf{Definitions} (m, V, W, H, N, S, K, R)

 $ m=4k $

 $ V = \{1,2,...,m\} $

 $ W = \{ \{a,b\} : a,b \in V \land a \neq b \} $

 $ H= \{ A \subset V : |A|=2k \} $

 $ \forall A \in H : N(A) = V \setminus A $

 $ \forall A \in H : S(A) = \{ \{a,b\} : a \in A \land b \in N(A) \}$

 $ K = \{ S(A) : A \in H \} $

 $ \forall T \subseteq K : R(T) = \{A \in H: S(A) \in T \} $

\textbf{Conclusions}

 $ |W|={{m}\choose{2}} $

 $ |H|={{m}\choose{2k}} $

 $ \forall X \in K : |X|=(2k)^2 $

 $ \forall A \in H : N(A) \in H $

 $ \forall A, B \in H : |A \cap B| = k \Longleftrightarrow |A \cap N(B)| = k $

 * Because $ |A \cap B| + |A \cap N(B)| = |A| = 2k $

 $ \forall A \in H: S(N(A)) = S(A) $

 $ \forall A, B \in H: S(A) = S(B) \Longrightarrow A=B \lor A=N(B) $

 $ \forall T \subseteq K : |R(T)| = 2 \cdot |T| $

 $ |K|=\frac{1}{2} \cdot |H| = \frac{1}{2} {{m}\choose{2k}} $

\textbf{Definitions} (L, d)

 $ L = \{ X \subseteq W \} $

 $ d = |W| - 1 $

\textbf{Conclusions}

 $ K \subset L $

 $ d = {{m}\choose{2}} - 1 $

 $\dim(L)=|W|$

 * Treating the sets in $L$ as \{0,1\}-vectors with index set $W$.

 $\dim(K) \leq \dim(L)-1 = d $

 * $K$ is in a hyperplane of $L$ because all $X \in K$ do have the same size.

\textbf{Definitions} ($\triangle$, dist, diam, M)

 $ \forall X, Y \in L : X \triangle Y = (X \setminus Y) \cup (Y \setminus X) $

 * The less-known general set operation \emph{symmetric difference}.

 $ \forall X, Y \in L : \mathrm{dist} (X,Y)=\sqrt{|X \triangle Y|}$

 * Euclidean distance for \{0,1\}-vectors represented by subsets.

 $ \forall T \subset K : \mathrm{diam}(T)
 = \max (\{ \mathrm{dist}(X,Y) : X,Y \in T \} \cup \{ 0 \} ) $

 $ M = \mathrm{diam}(K)$

\textbf{Conclusions}

 $\exists A,B \in H : |A \cap B| = k $

 * For instance  $A=\{1,2,...,2k\} \land B= \{k+1,k+2,...,3k\} $

 $ \forall A, B \in H : \mathrm{dist}(S(A),S(B)) = M \Longrightarrow |A \cap B| = k $

 * Not obvious, not proved in \cite{K}, but we give a proof at the end of this subsection.

\textbf{Definition} (Z)

 $ Z= \{ F \subseteq H : (\forall A, B \in F : |A\cap B| \neq k) \}$

\textbf{Conclusions}

 $\forall T \subset K : \mathrm{diam}(T) < M \Longrightarrow R(T) \in Z $

 $\forall F \in Z : |F| \leq 2 \cdot {{m-1}\choose{k-1}} $

 * Instantiation of the Frankl-Wilson theorem.

 $\forall T \subseteq K : \mathrm{diam}(T) < M \Longrightarrow |R(T)| \leq 2 \cdot {{m-1}\choose{k-1}}
  \land |T| = \frac{|R(T)|}{2} \leq {{m-1}\choose{k-1}} $

\textbf{Definitions} (G, q)

$ G = \{ P \subseteq \{ T \subseteq K : \mathrm{diam}(T) < M \} : \bigcup_{T \in P} T = K \} $

$ q= \frac{1}{2} {{m}\choose {2k}} / {{m-1}\choose {k-1}} $

\textbf{Conclusion}

$\forall P \in G : |P| \geq |K| / {{m-1}\choose{k-1}} = q $

\subsubsection{Textual summary of the result}

The point set $K$ lies in the Euclidean space of dimension
 $d={{m}\choose {2}}-1$ and can not be covered by less than

$$q=\frac{\frac{1}{2} {{m}\choose {2k}}}{{{m-1}\choose {k-1}}}
 = \frac{ {{m}\choose {2k}}}{ 2\cdot { {m-1}\choose {k-1}}} $$

subsets of smaller diameter. If $q > d+1 $ then $K$ is a counterexample to
Borsuk's conjecture, in dimension $d$ and all higher dimensions below $q-1$.
Otherwise we can not prove it this way.

Compared to the respective formula given in \cite{K}, our formula for the
lower bound delivers doubled values. This results in extended possibilities
with respect to direct verifications of constructed point sets as counterexamples.

\subsubsection{Modified formulas and bounds for q}

Observe that for all positive integers $a$ and $b$, $a \geq b$, it is

$$
 \frac{a}{b} \cdot {{a-1}\choose {b-1}} =
 \frac{a}{b} \cdot \frac{(a-1)!}{(b-1)!\cdot((a-1)-(b-1))!} =
 \frac{a\cdot(a-1)!}{b\cdot (b-1)!\cdot(a-b)!} =
 \frac{a!}{b!\cdot(a-b)!} = {{a}\choose {b}} $$

\medskip

Here we set $a=m=4k$, $b=k$ and get
 $ 4 \cdot {{m-1}\choose {k-1}} = {{m}\choose {k}} $ and can derive

$$ q/2 =\frac{ {{m}\choose {2k}} }{{{m}\choose {k}}} =
  \frac{ \frac{m!}{(2k)!\cdot(m-2k)!} }{ \frac{m!}{k!\cdot(m-k)!} } =
  \frac{ {k!\cdot(m-k)!}  }{ {(2k)!\cdot(m-2k)!} } =
  \frac{ {k!\cdot(3k)!}  }{ {(2k)!\cdot(2k)!} } =
  \frac{ \frac{ (3k)!}{(2k)!} }{ \frac{(2k)!}{k!} } =
  \frac{ \frac{ (3k)!}{k! \cdot(2k)!} }{ \frac{(2k)!}{k! \cdot k!} } =
  \frac{ {{3k}\choose {k}} }{{{2k}\choose {k}}}
  $$

We can use the last expression for $q/2$ this way
$$ q = 2 \cdot \frac{ {{3k}\choose {k}} }{{{2k}\choose {k}}}
   = \frac{ {{3k}\choose {k}} }{{{2k-1}\choose {k-1}}}
$$
but also interpret the third-last one this way
$$ q
= 2 \cdot \frac{ {\prod_{i=2k+1}^{3k} i} }{ \prod_{i=k+1}^{2k} i }
= \frac{ {\prod_{i=2k+1}^{3k} i} }{ \prod_{i=k}^{2k-1} i }
$$
Both expressions include relatively few multiplications and just one division.

\medskip

In order to estimate bounds, we use the first expression and derive

$$ q
= 2 \prod_{i=k+1}^{2k} \frac{k+i}{i}
= 2 \prod_{i=k+1}^{2k} (1+ \frac{k}{i})
$$

The second expression proves that the iterated factor gets its minimum
and its maximum iff $i=2k$ and $i=k+1$, respectively, but using the iterated
factor of the first expression shortens the derivation.

 $$ 2k \mapsto \frac{k+2k}{2k}=\frac{3}{2} \ \ , \ \
 k+1 \mapsto \frac{k+k+1}{k+1}=\frac{2k+2-1}{k+1}=2-\frac{1}{k+1}<2$$

Consequently $$ 2 \cdot (\frac{3}{2})^k \leq q \leq 2 \cdot (2-\frac{1}{k+1})^k<2^{k+1} $$

With respect to $k$, $d$ grows quadratically and $q$ grows exponentially.

\subsubsection{Actual calculation results and covered dimensions}

$ k=11, d=945, q \approx 548.70 $

* Point set not verified as counterexample

$ k=13,  d=1325, q \approx 1561.91 $

* Counterexample for $1325 \leq d \leq 1560$

$ k=16, d=2015, q \approx 7502.65 $

* Counterexample for $2015 \leq d \leq 7501$

$ k=17, d=2277, q \approx 12659.44 $

* Counterexample for $2277 \leq d \leq 12658 $

$ k=29, d=6669, q \approx 6745998.54 $

$ k=31, d=7625, q \approx 19209098.12 $

$ k=32, d=8127, q \approx 32414445.61 $

\medskip

To provide counterexamples for any $d \geq 2015$ it suffices
to consider for $k$ the powers of two from 16 onwards and one of the prime
powers from 17 to 29 (to avoid a gap before $k=32$ applies).

\subsubsection{Proof that the maximality of $\mathrm{dist}(S(A),S(B))$ implies $|A\cap B|=k$ }

This proof is aimed to justify the above (3.2.2) conclusion

$ \forall A, B \in H : \mathrm{dist}(S(A),S(B)) = M \Longrightarrow |A \cap B| = k $

within its environment (preceding propositions), i.e. the premisses and
definitions stated herein are additional.

\textbf{Premisses} (A, B)

$A \in H$

$B \in H$

\textbf{Definitions} (a, b, g, p)

$a = S(A)$

$b = S(B)$

$g = |a \triangle b|$

$p = |A\cap B|$

\textbf{Conclusions}

Because $\mathrm{dist}(a,b)=\sqrt{g} $, it suffices to prove that $g$
reaches its maximum if and only if $p=k$ .

$ g = |(a \setminus b) \cup (b \setminus a)| $

$ (a \setminus b) \cap (b \setminus a) = \emptyset $

$ g = |a \setminus b| + |b \setminus a|
 = |a| - |a \cap b| + |b| - |b \cap a| = 2 \cdot (2k)^2 - 2\cdot |a \cap b| $

Now we know that $g$ is maximal if and only if $|a \cap b|$ is minimal and
it just remains to show that $|a \cap b|$ reaches its minimum if and
only if $p=k$.

\textbf{Conclusions}

$ a = \{ \{x,y\} : x \in A \land y \in N(A) \} $

$ b = \{ \{x,y\} : x \in B \land y \in N(B) \} $

$ a \cap b = \{ \{x,y\} :
 x \in A \cap B \land y \in N(A) \cap N(B) \lor
 x \in A \cap N(B) \land y \in N(A) \cap B
 \} $

$ a \cap b =
 \{ \{x,y\} : x \in A \cap B \land y \in N(A) \cap N(B) \} \cup
 \{ \{x,y\} : x \in A \cap N(B) \land y \in N(A) \cap B \} $

$ |a \cap b| =
 |\{ \{x,y\} : x \in A \cap B \land y \in N(A) \cap N(B) \}| +
 |\{ \{x,y\} : x \in A \cap N(B) \land y \in N(A) \cap B \}| $

$ |a \cap b| = |A\cap B| \cdot |N(A)\cap N(B)| + |A\cap N(B)| \cdot |N(A)\cap B| $

$|A\cap N(B)| = |A| - |A \cap B| = 2k - p$

$|N(A)\cap B| = |B| - |B \cap A| = 2k - p$

$|N(A) \cap N(B)| = |N(A)| - |N(A) \cap B| = 2k - (2k - p) = p $

$ |a \cap b| = p^2 + (2k-p)^2 = p^2 + 4k^2 - 4kp + p^2
 = 2p^2 + 4k^2 - 4kp$

\textbf{Definition} (v)

$ v=p-k $

\textbf{Conclusions}

$ p=k+v $

$|a \cap b| = 2(k+v)^2 + 4k^2 - 4k(k+v)
= 2k^2+4kv+2v^2 + 4k^2 - 4k^2 - 4kv = 2k^2+2v^2 $

Thus, the minimum value of $|a \cap b|$ is $2k^2$, and it is reached if and
only if $v=0$, i.e. $p=k$.

\vspace{0.1in}

Author's eMail address: thomas.jenrich@gmx.de


\begin{thebibliography}{00}

\bibitem{B}
K. Borsuk, \emph{Drei S\"{a}tze \"{u}ber die $n$-dimensionale euklidische
Sph\"{a}re}, Fund. Math., 20 (1933), 177-190.

\bibitem{F}
P. Frankl and R. Wilson,
 \emph{Intersection theorems with geometric consequences},
 Combinatorica 1 (1981), 357-368.

\bibitem{K}
J. Kahn and G. Kalai, \emph{A counterexample to Borsuk's conjecture},
Bull. Amer. Math. Soc. 29 (1993), 60-62.
arXiv:math.MG/9307229v1

\bibitem{L}
D. Larman, \emph{Open problem 6}, Convexity and Graph Theory (M. Rozenfeld and
J. Zaks, eds.), Ann. Discrete Math., vol. 20, North-Holland,
Amsterdam and New York, 1984, p. 336.

\bibitem{W}
B. Wei{\ss}bach, \emph{Sets with large Borsuk number},
 Beitr\"age Algebra Geom. 41 (2000), 417-423.

\end{thebibliography}
\end{document}